\input amstex
\documentstyle{amsppt}
\magnification=\magstep1
\hcorrection{0in}
\vcorrection{0in}
\pagewidth{6.3truein}
\pageheight{9truein}

\topmatter
\title Decomposition of log crepant birational morphisms between log terminal surfaces
\endtitle
\author Shigetaka Fukuda
\leftheadtext{S. FUKUDA}
\rightheadtext\nofrills{\rm log crepant birational morphisms}
\endauthor
\address
Faculty of Education,
Gifu Shotoku Gakuen University,
Yanaizu-cho,
hashima-gun,
Gifu prefecture,
501-6194 Japan
\endaddress
\email fukuda\@ha.shotoku.ac.jp
\endemail
\abstract
We prove that every log crepant birational morphism between log terminal surfaces is 
decomposed into log-flopping type divisorial contraction morphisms and log blow-downs. 
Repeating these two kinds of contractions we reach a minimal log minimal surface from any log 
minimal surface.
\endabstract
\endtopmatter

\document
\head 1. Introduction
\endhead
All varieties are defined over the complex number field $\Bbb C$ in this paper.
We generally use the notation and terminology of [1].

\definition{Definition 1}
Let $X$ be a normal algebraic surface and $D$ an effective $\Bbb Q$-divisor on $X$ such that 
$\lceil D \rceil$ is reduced.
A log surface $(X,D)$ is said to be {\it log terminal} if the following
conditions are satisfied:

(1) $K_X+D$ is $\Bbb Q$ -Cartier.

(2) There exists a log resolution $f:Y \to X$ such that $K_Y+f^{-1}_*D = f^*(K_X+D)+ \sum 
a_jE_j$ for $a_j \in \Bbb Q$ with $a_j > -1$.

A proper birational morphism $h:(\tilde{X},\tilde{D}) \to (\bar{X},\bar{D})$ between log 
canonical surfaces is said to be a {\it log crepant} birational morphism if $K_{\tilde{X}} + 
\tilde{D} = h^*(K_{\bar{X}}+\bar{D})$.
\enddefinition

\remark{Remark}
The following facts are well known to experts:

(1) The notions of log terminal, divisorial log terminal and weakly Kawamata log terminal are 
equivalent in case of surfaces.

(2) Every complete log terminal surface is $\Bbb Q$-factorial and projective.
Thus every proper surjective morphism from any log terminal surface is projective.
\endremark

\definition{Definition 2}
(Minimal log minimal surface)\quad
Let $(\tilde{X},\tilde{D})$ be a log terminal surface and $g:\tilde{X} \to T$ a proper 
morphism onto a variety $T$.
In this case, $(\tilde{X},\tilde{D})$ is denoted by $(\tilde{X},\tilde{D})/T$ and called a log 
terminal surface$/T$. 
We call $g$ the {\it structure morphism} of $\tilde{X} /T$.

If $K_{\tilde{X}} + \tilde{D}$ is $g$-nef, we say $(\tilde{X},\tilde{D})/T$ is a {\it log 
minimal} surface$/T$ (When $T=$ Spec $\Bbb C$, we simply call $(\tilde{X},\tilde{D})$ a {\it 
log minimal} surface.).

Let $(\bar{X},\bar{D})$ be another log terminal surface$/T$.
A log crepant birational morphism $h:(\tilde{X},\tilde{D}) \to (\bar{X},\bar{D})$ being 
compatible with the structure morphisms of $\tilde{X} /T$ and $\bar{X} /T$ is denoted by 
$h:(\tilde{X},\tilde{D}) /T \to (\bar{X},\bar{D}) /T$ and called a {\it log crepant} 
birational morphism$/T$.

A log minimal surface $(\tilde{X},\tilde{D}) /T$ is said to be a {\it minimal log minimal} 
surface$/T$ if every log crepant birational morphism$/T$ from $(\tilde{X},\tilde{D}) /T$ is an 
isomorphism (When $T=$ Spec $\Bbb C$, we simply call $(\tilde{X},\tilde{D})$ a {\it minimal 
log minimal} surface.).
\enddefinition

\definition{Construction}
(Log-flopping type divisorial contraction)\quad
Let $(\tilde{X},\tilde{D}) /T$ be a log minimal surface$/T$ with structure morphism $g$.
Assume that a $g$-exceptional curve $C$, with $(K_{\tilde{X}} + \tilde{D}, C) = 0$ and $C^2 < 
0$, does not contain the center of $\nu$ on $\tilde{X}$ for any divisor $\nu$ of 
Rat$(\tilde{X})$ with discrepancy $-1$ with rspect to $(\tilde{X},\tilde{D})$ (We say $C$ is a 
{\it log-flopping type} divisor with respect to  $(\tilde{X},\tilde{D}) /T$.).
Then $(\tilde{X},\tilde{D}+\epsilon C)$ is log terminal for any sufficiently small positive 
rational number $\epsilon$.
Thus $C$ spans an extremal ray $R_C$ for $(\tilde{X},\tilde{D}+\epsilon C) /T$.
Consequently we have the divisorial contraction morphism $h:\tilde{X} /T \to \bar{X} /T$ of 
the extremal ray $R_C$.
Putting $\bar{D}:=h_* \tilde{D}$, we have a log crepant birational morphism 
$h:(\tilde{X},\tilde{D}) /T \to (\bar{X},\bar{D}) /T$ between log minimal surfaces$/T$ (The 
morphism $h$ is said to be a {\it log-flopping type} divisorial contraction morphism$/T$. When 
$T=$ Spec $\Bbb C$, it is simply called a {\it log-flopping type} divisorial contraction 
morphism.).
\enddefinition

\definition{Definition 3}
Let $(\tilde{X},\tilde{D}) /T$ be a log minimal surface$/T$.
If there is no log-flopping type divisor with respect to  $(\tilde{X},\tilde{D}) /T$, we say 
$(\tilde{X},\tilde{D}) /T$ is a {\it log-flopping-type-divisors-contracting-process minimal} 
surface$/T$.
\enddefinition

\definition{Definition 4}
(Log blow-down)\quad
A log crepant birational morphism $h:(\tilde{X},\tilde{D}) \to (\bar{X},\bar{D})$ between log 
terminal surfaces is called the {\it log blow-down} to a closed point $x \in \bar{X}$ if there 
is an open neighborhood $U$ of $x$ with the following properties:

(1)$U$ is smooth.

(2)$\bar{D} \vert_U = \Delta_1 +\Delta_2$ that is a reduced simple normal crossing divisor 
with $\Delta_1 \cap \Delta_2 = \{ x \}$.

(3)$h^{-1}$ is the blow-up at $x$.

(4)$h \vert_{h^{-1}(\bar{X} \setminus \{ x \})}$ is an isomorphism.
\enddefinition

Now we state our main theorems.

\proclaim{Main Theorem I}
{\rm (Decomposition)}
Let $\varphi : (X_1,D_1) \to (X_2,D_2)$ be a log crepant birational morphism between log 
terminal surfaces.
Starting with $(X_1,D_1) /X_2$ , after a sequence of log-flopping type divisorial 
contractions$/X_2$ ,we end up with a log-flopping-type-divisors-contracting-process minimal 
surface $(X_{fm},D_{fm}) /X_2$.
Furthermore the structure morphism of $(X_{fm},D_{fm}) /X_2$ is a log crepant birational 
morphism from $(X_{fm},D_{fm})$ to $(X_2,D_2)$ and it is a composite of log blow-downs.
\endproclaim

We prove the theorem above in the next section.

\proclaim{Main Theorem II}
{\rm (Arrival at a minimal log minimal surface)}\quad
Let $(X,D)$ be a log minimal surface.
Then after repetitions of log-flopping type divisorial contractions and log blowing-downs, we 
reach a minimal log minimal surface.
\endproclaim

\demo{Proof}
We note that the Picard number strictly decreases after a log-flopping type divisorial 
contraction and also after a log blowing-down.
As well, in Definition 3, if $\tilde{X}$ is complete and $K_{\tilde{X}} + \tilde{D}$ is nef, 
then every log-flopping type divisorial contraction$/T$ is a log-flopping type divisorial 
contraction (/Spec $\Bbb C$).
Thus the assertion follows from Main Theorem I.
\qed
\enddemo

\head 2. Proof of Theorem I
\endhead

\proclaim{Lemma 1}
{\rm ([2, 1.1.])}
Assume that $(X,D)$ is a log terminal surface and that $x$ is a closed point on $X$ such that 
$\{ x \}$ is the center {\rm center}$_X (\nu)$ on $X$ of a divisor $\nu$ of {\rm Rat}$(X)$ 
with discrepancy $-1$ with respect to $(X,D)$.
Then there exists an open neighborhood $U$ of $x$ such that $U$ is smooth and $D \vert_U = 
\Delta_1 +\Delta_2$ which is a reduced simple normal crossing divisor with $\Delta_1 \cap 
\Delta_2 = \{ x \}$.
\endproclaim

\proclaim{Lemma 2}
Let $h: \tilde{X} \to \bar{X}$ be a proper birational morphism between smooth surfaces and 
$\bar{D}$ a reduced simple normal crossing divisor on $\bar{X}$.
Then $h^{-1}_* \bar{D} + \sum \{E \vert E$ {\rm is an} $h$-{\rm exceptional prime divisor} $\}
$ is a reduced simple normal crossing divisor.
\endproclaim

\proclaim{Proposition 1}
Let $h:(\tilde{X},\tilde{D}) \to (\bar{X},\bar{D})$ be a log crepant birational morphism 
between log terminal surfaces and let $C$ be an $h$-exceptional curve such that $C$ is not an 
irreducible component of
$\lfloor \tilde{D} \rfloor$.
Then $C$ does not contain {\rm center}$_{\tilde{X}} (\nu)$ for any divisor $\nu$ of {\rm 
Rat}$(\tilde{X})$ with discrepancy $-1$ with respect to $(\tilde{X},\tilde{D})$.
\endproclaim

\demo{Proof}
We will derive a contradiction assuming that, for some divisor $\nu$ of {\rm Rat}$(\tilde{X})$ 
with discrepancy $-1$ with respect to $(\tilde{X},\tilde{D})$, $C$ contains {\rm 
center}$_{\tilde{X}} (\nu)$.
Then {\rm center}$_{\tilde{X}} (\nu) = \{ p \}$ for some closed point $p$ on $\tilde{X}$.
Hence, from Lemma 1, there exists an open neighborhood $U$ of $p$ such that $U$ is smooth and 
$\tilde{D} \vert_U = \tilde{\Delta}_1 + \tilde{\Delta}_2$ which is a reduced simple normal 
crossing divisor with $\tilde{\Delta}_1 \cap \tilde{\Delta}_2 = \{ p \}$.
We note that $\{ h(p) \} =$ {\rm center}$_{\bar{X}}(\nu)$.
Hence from Lemma 1 again, there exists an open neighborhood $V$ of $h(p)$ such that $V$ is 
smooth and $\bar{D} \vert_V = \bar{\Delta}_3 + \bar{\Delta}_4$ which is a reduced simple 
normal crossing divisor with $\bar{\Delta}_3 \cap \bar{\Delta}_4 = \{ h(p) \}$.
Here $(\tilde{\Delta}_1 + \tilde{\Delta}_2 + C) \vert_{U \cap h^{-1} (V)} \leq ((h 
\vert_{h^{-1} (V)})^{-1}_* (\bar{\Delta}_3 + \bar{\Delta}_4) + \sum \{E \vert E$ {\rm is an} 
$h$-{\rm exceptional prime divisor} $\}) \vert_{U \cap h^{-1} (V)}$.
But $(\tilde{\Delta}_1 + \tilde{\Delta}_2 + C) \vert_{U \cap h^{-1} (V)}$ is not a reduced 
simple normal crossing divisor.
This is a contradiction, by Lemma 2.
\enddemo

\proclaim{Proposition 2}
Let $g:(X,D) \to (\bar{X},\bar{D})$ be a log crepant birational morphism between log terminal 
surfaces such that every $g$-exceptional
prime divisor is an irreducible component of $\lfloor D \rfloor$.
Then $g$ is a composite of log blow-downs.
\endproclaim

\demo{Proof}
Let $p_i$ $(i \in I)$ be the closed points on $\bar{X}$ such that $p_i$ is the generic point 
of {\rm center}$_{\bar{X}} (\nu)$ for some divisor $\nu$ of {\rm Rat}$(\bar{X})$ with 
discrepancy $-1$ with respect to $(\bar{X},\bar{D})$.
We note that, over $\bar{X} \setminus B$, $g$ is an isomorphism where $B = \{ p_i \vert \quad 
i \in I \}$.
We take a log resolution $f:Y \to X$ of $(X,\Delta)$ as in Definition 1.
Then every divisor on $Y$ with discrepancy $-1$ with respect to $(\bar{X},\bar{D})$ is an 
irreducible component of $f^{-1}_* \lfloor D \rfloor$.
Now we consider the morphism $\psi:=gf:Y \to \bar{X}$.
Then from Lemma 1 there exists an open neighborhood $U_i$ of $p_i$ with the following 
properties:

(1) $U_i$ is smooth.

(2) $\bar{D} \vert_{U_i} = \Delta_{i_1} +\Delta_{i_2}$ that is a reduced simple normal 
crossing divisor.

(3) $p_i$ is the generic point of $\Delta_{i_1} \cap \Delta_{i_2}$.

(4) $\psi \vert_{\psi ^{-1} \left( U_i \right)}$ is a proper birational morphism between 
amooth surfaces and over $U_i \setminus B$ it is an isomorphism.

Now, with $\psi ^{-1} \left( U_i \right)$ toward $U_i$, we start the process of contracting 
$(-1)$-curves being exceptional over $U_i$ with discrepancies $>-1$ with respect to $(U_i, 
\bar{D} \vert_{U_i})$.
Then we end up with a smooth surface $V_i$ such that every $(-1)$-curve being exceptional over 
$U_i$ is with discrepancy $-1$ with respect to $(U_i, \bar{D} \vert_{U_i})$.

Next with $V_i$ toward $U_i$, we start the process of contracting $(-1)$-curves being 
exceptional over $U_i$ with discrepancies $-1$ with respect to $(U_i, \bar{D} \vert_{U_i})$.
At every stage of this process a $(-1)$-curve being exceptional over $U_i$ with discrepancy 
$-1$ with respect to $(U_i, \bar{D} \vert_{U_i})$ contracts to a point that is the 
intersection of two prime divisors with discrepancies $-1$ with respect to $(U_i, \bar{D} 
\vert_{U_i})$.
From Lemma 2, there exists no other prime divisor that is exceptional over $U_i$ and passes 
through this point.
Thus, during this process, a $(-1)$-curve being exceptional over $U_i$ with discrepancy $>-1$ 
with respect to $(U_i, \bar{D} \vert_{U_i})$ is not born.
Therefore after this process we reach $U_i$.

Consequently every curve on $V_i$ that is exceptional over $U_i$ is with discrepancy $-1$ with 
respect to $(U_i, \bar{D} \vert_{U_i})$.
Here we note that the $V_i$ and $\bar{X} \setminus B$ patch together to a complete surface 
$M$.
By the argument above and the choice of a log resolution $f$, $M \setminus P_M$ is isomorphic 
to $X \setminus P_X$ where $P_M$ (resp\. $P_X$) is a closed set composed of a finite number of 
closed points on $M$ (resp\. $X$).
As a result $M$ is isomorphic to X, from Zariski's Main Theorem.
\enddemo

\definition{Proof of Main Theorem I}
Proposition 1, Construction and Proposition 2 imply the assertion.
\enddefinition

\enddocument